\documentclass[amsmath,amssymb,pra]{revtex4-1}
\usepackage{fullpage}
\pdfoutput=1
\usepackage{amsfonts}
\usepackage{amssymb}
\usepackage{amsmath}
\usepackage{graphicx}
\usepackage{bbm}
\usepackage{overpic}
\usepackage{subfigure}
\usepackage{xcolor}

\newtheorem{lemma}{Lemma}
\newtheorem{remark}{Remark}
\newtheorem{example}{Example}
\newtheorem{definition}{Definition}
\newtheorem{proposition}{Proposition}

\newtheorem{property}{Property}

\newcommand{\proof}{\noindent {\it Proof.\ }}
\newcommand{\scalar}[2]{\langle #1 | #2 \rangle}
\newcommand{\ketbra}[2]{| #1 \rangle \langle #2 |}
\newcommand{\ket}[1]{| #1 \rangle}
\newcommand{\bra}[1]{\langle #1 |}

\newcommand{\kron}{\otimes}
\newcommand{\tr}{\mathrm{tr}}
\newcommand{\1}{{\rm 1\hspace{-0.9mm}l}}
\newcommand{\Id}{\1}
\newcommand{\halmos}{\hfill $\blacksquare$\newline}
\newcommand{\R}{\ensuremath{\mathbb{R}}}
\newcommand{\setR}{\R}
\newcommand{\diag}{\mathrm{diag}}

\newcommand{\Cplx}{\mathbb{C}}

\newcommand{\MM}{\mathbb{M}}

\newcommand{\LambdaRm}{\mathrm{\Lambda}}
\newcommand{\NumRange}[1]{\ensuremath{\LambdaRm({#1})}}

\newcommand{\ProductNumRange}[1]{\ensuremath{\LambdaRm^{\!\otimes}\! \left( #1 \right)}}

\newcommand{\lambdaProdMin}{\lambda_{\rm min}^{\otimes}}
\newcommand{\lambdaProdMax}{\lambda_{\rm max}^{\otimes}}
\newcommand{\ie}{\emph{i.e.}}
\newcommand{\eg}{\emph{e.g.}}
\newcommand{\etal}{\emph{et al.}}
\newcommand{\mplus}{\boxplus}
\newcommand{\mprod}{\boxtimes}

\renewcommand{\bar}[1]{{#1}^*}
\usepackage{xcolor}

\begin{document}

\title{Product numerical range in a space with tensor product structure}
\author{Zbigniew Pucha{\l}a}
\email{z.puchala@iitis.gliwice.pl}
\author{Piotr Gawron}
\author{Jaros{\l}aw Adam Miszczak}
\affiliation{Institute of Theoretical and Applied Informatics, Polish Academy
of Sciences, Ba{\l}tycka 5, 44-100 Gliwice, Poland}
\author{{\L}ukasz~Skowronek}
\affiliation{Instytut Fizyki im. Smoluchowskiego, Uniwersytet
Jagiello{\'n}ski, Reymonta 4, 30-059 Krak{\'o}w, Poland }
\author{Man-Duen~Choi}
\affiliation{Department of Mathematics, University of Toronto, 
Toronto, Ontario, Canada M5S 2E4}
\author{Karol \.Zyczkowski}
\affiliation {Instytut Fizyki im. Smoluchowskiego, Uniwersytet
Jagiello{\'n}ski,
ul. Reymonta 4, 30-059 Krak{\'o}w, Poland}
\affiliation{Centrum Fizyki Teoretycznej, Polska Akademia Nauk, Al.
Lotnik{\'o}w
32/44, 02-668 Warszawa, Poland}

\begin{abstract}
We study operators acting on a~tensor product Hilbert space and investigate
their product numerical range, product numerical radius and separable numerical
range. Concrete bounds for the product numerical range for Hermitian operators
are derived. Product numerical range of a~non-Hermitian operator forms a~subset
of the standard numerical range containing the barycenter of the spectrum. While the latter set is convex, the product
range needs not to be convex
nor simply connected. The product numerical range
of a~tensor product is equal to the Minkowski product of numerical ranges of
individual factors.
\end{abstract}

\date{August 20, 2010}

\maketitle

%

\section{Introduction}

Let $X$ be an operator acting on an $N$-dimensional Hilbert space ${\cal H}_N$.
Let $\NumRange{X}$ denote its \emph{numerical range}, \ie\ the set of all
$\lambda$ such that there exists a~normalized state $\ket{\psi}\in
\mathcal{H}_N$, $|| \psi || = 1$, which satisfies $ \bra{\psi} X \ket{\psi} =
\lambda$.

In this work we study an analogous notion defined for operators
acting on a composite Hilbert space with a tensor product structure.
Consider first a bi--partite Hilbert space, 
\begin{equation}
\mathcal{H}_N = \mathcal{H}_K \otimes \mathcal{H}_M ,
\label{HKM}
\end{equation}
of a composite dimension $N=KM$. 

\begin{definition}[Product numerical range]
Let $X$ be an operator acting on the composite Hilbert space (\ref{HKM}).
We define the \emph{product numerical range} \ProductNumRange{X} 
of~$X$, with respect to the tensor product structure of $\mathcal{H}_N$, as
\begin{equation}\label{lnr}
  \ProductNumRange{X} = \left\{ \bra{\psi_A \otimes \psi_B} X \ket{\psi_A\otimes\psi_B} 
: \ket{\psi_A}\in \mathcal{H}_K, \ket{\psi_B}\in \mathcal{H}_M \right\},
\end{equation}
where
$\ket{\psi_A} \in {\cal H}_K$  and $\ket{\psi_B} \in {\cal H}_M$ are normalized. 
\end{definition}

\begin{definition}[Product numerical radius]
Let $\mathcal{H}_N = \mathcal{H}_K \otimes \mathcal{H}_M$ be a~tensor product 
Hilbert space. We define the \emph{product numerical radius} $r^{\otimes}(X)$
of~$X$, with respect to this tensor product structure, as
\begin{equation}
r^{\otimes}(X) = \max\{|z|:z\in\ProductNumRange{X}\}.
\end{equation}
\end{definition}

The notion of numerical range of a~given operator, also called ``field of
values'' \cite[Chapter 1]{HJ2}, has been extensively studied during the last few
decades \cite{gustav,ando94numerical,li95cnumerical}
and its usefulness in quantum theory has been emphasized \cite{KPLRS09}.
Several generalizations of numerical range  are known -- 
see \eg\ \cite[Section 1.8]{HJ2}. In particular, Marcus introduced the notion 
of \emph{decomposable numerical range} \cite{Ma73,MW80}, the properties of which 
are a~subject of considerable interest
\cite{BLP91,LZ01}.

The product numerical range, which forms the central point of this work, can be
considered as a~particular case of the decomposable numerical range defined for
operators acting on a~tensor product Hilbert space. This notion may also be
considered as a~numerical range \emph{relative} to the proper
subgroup $U(K)\times U(M)$ of the full unitary group $U(KM)$.

In papers \cite{dirr08relative,DFY08,thomas08significance,SHGDH08} the 
same
object was called \emph{local numerical range} in view of notation common in
quantum mechanics. This name seems to be natural for the physicists audience,
but to be more consistent with the mathematical terminology we will use in this
paper the name \emph{product numerical range}, although a~longer version ``local
product numerical range'' would be even more accurate.

In a~recent paper of Dirr \etal\ \cite{dirr08relative} some geometric
properties of the product numerical range and product $C$-numerical range were
investigated. Another paper of the same group \cite{thomas08significance}
demonstrates the possible application of these concepts in the theory of quantum
information and the theory of quantum control. Product numerical range of
unitary operators was very recently used by Duan \etal{} to tackle the problem
of local distinguishability of unitary operators \cite{DFY08}.
Knowing  product numerical range of a Hermitian operator one could 
solve other important problems in the theory of quantum information,
as establishing whether a given quantum map is positive,
or obtaining bounds for the minimum output entropy 
of a~quantum channel \cite{Gxx10}.

The main goal of this paper is to stimulate research on product numerical ranges.
We derive several bounds for the product numerical range of a Hermitian operator
defined on a space with a two--fold tensor structure, 
which corresponds to a bi-partite physical system.
In the non-Hermitian case we use the relation between 
the product numerical range of a tensor product and the Minkowski product
of numerical ranges to establish a general bound
for the product numerical range 
based on the operator Schmidt decomposition.
We show that a tensor product of two operators, acting on a Hilbert space
with a two--fold tensor structure, has a simply connected product numerical range.
A similar property does not hold for operators acting on a space with a
larger number of factors. We introduce a class of 
product diagonalizable operators, for which a convenient method to parameterize
product numerical range is proposed.

Although this work leaves several problems related to product numerical 
range unsolved, we believe it could point out directions
for further mathematical research,
which will find direct applications in the theory of quantum information.
In order to invite reader to contribute to this field we 
conclude the paper with a list of exemplary open problems.

\section{Properties of product numerical range}
\label{sec:nonHermitian}
In this section we are going to consider arbitrary operators acting on a
bipartite Hilbert space~(\ref{HKM}). 

\subsection{General case}
It is not difficult to establish the basic properties of the product numerical
range which are independent of the partition of the Hilbert space and of the
structure of the operator. We list them below leaving some simple items without
a~proof.
\subsubsection{Basic properties}
We begin this section with some simple topological facts concerning product
numerical range for general operators.

\begin{property}
	Product numerical range forms a~connected set in the complex plane.
\end{property}
\proof
The above is true because product numerical range is a~continuous image of 
a~connected set.
\halmos
\begin{property}[Subadditivity]\label{prop:subadditivity}
Product numerical range is subadditive. For all $A, B\in \MM_n$
	\begin{equation}
	\ProductNumRange{A+B}\subset \ProductNumRange{A}+\ProductNumRange{B}.
	\end{equation}
\end{property}

\begin{property}[Translation]
	For all $A\in \MM_n$ and $\alpha\in \Cplx$
	\begin{equation}
	\ProductNumRange{A+\alpha \Id}=\ProductNumRange{A}+\alpha.
	\end{equation}
\end{property}

\begin{property}[Scalar multiplication]
	For all $A\in \MM_n$ and $\alpha\in \Cplx$
	\begin{equation}
	\ProductNumRange{\alpha A}=\alpha\ProductNumRange{A}.
	\end{equation}
\end{property}

\begin{property}[Product unitary invariance]
	For all $A\in \MM_{mn}$ 
	\begin{equation}
	\ProductNumRange{(U\otimes V)A(U\otimes V)^\dagger}=\ProductNumRange{A},
	\end{equation}
for unitary $U\in \MM_m$ and $V\in \MM_n$.
\end{property}

\begin{property}
Let $A \in \MM_m$ and $B\in \MM_n$ 
\begin{enumerate}
\item If one of them is normal then the numerical range 
of their tensor product coincides with the convex hull of the
product numerical range,
\begin{equation}
\NumRange{A \kron B} = \mathrm{Co}(\ProductNumRange{A \kron B}).
\end{equation}
\item If $e^{i \theta} A$ is positive semidefinite for some $\theta \in[0,2 \pi)$, then
\begin{equation}
\NumRange{A \kron B} = \ProductNumRange{A \kron B}.
\end{equation}
\end{enumerate}
\end{property}
\proof
This property can be proven using Lemma~\ref{lemma:KronToMin} stated in the
following subsection and Theorem 4.2.16 in~\cite{HJ2}.
\halmos

Let $H(A)=\frac{1}{2}(A+A^\dagger)$ and $S(A)=\frac{1}{2}(A-A^\dagger)$.
\begin{property}[Projection]
	For all $A\in \MM_{n}$ 
	\begin{equation}
	\ProductNumRange{H(A)}=\mathrm{Re}\ \ProductNumRange{A}
	\end{equation}
	and
	\begin{equation}
	\ProductNumRange{S(A)}=i\,\mathrm{Im}\ \ProductNumRange{A}.
	\end{equation}
\end{property}

\begin{property}
The product numerical range does not need to be convex. 
\end{property}
\proof 
Consider the following simple example.

\begin{figure}[ht!]
 \begin{center}
  \includegraphics[scale=0.6]{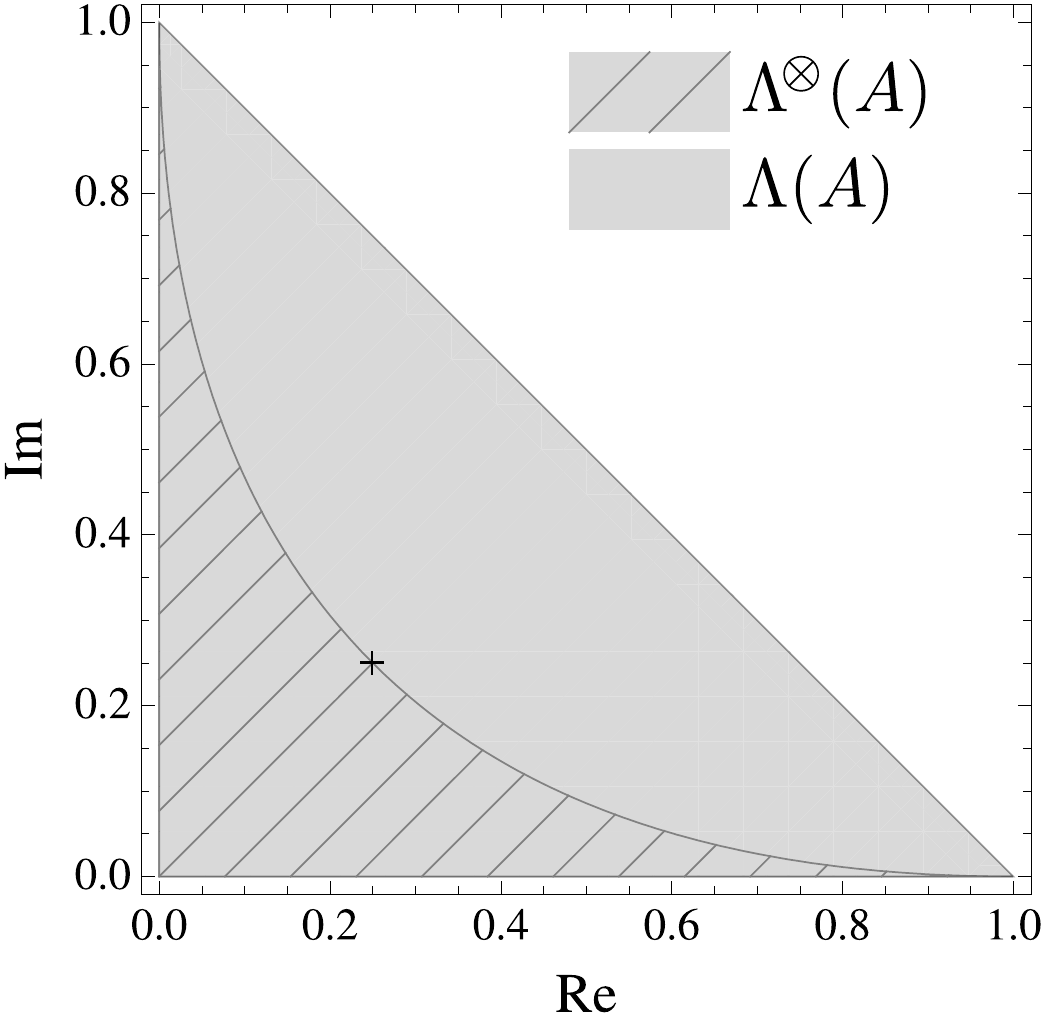}
  \caption{The comparison of the numerical range (gray triangle)
  and the product numerical range (dashed set) for 
  matrix $A$ defined in Eq. \eqref{eqn:example-nonconvex}.}
	\label{fig:example-nonconvex}
 \end{center}
\end{figure}

\begin{example}\label{example-nonconvex}
Let
\begin{equation}
\label{eqn:example-nonconvex}
A=
\left(
\begin{array}{cc}
1 & 0\\
0 & 0
\end{array}
\right)
\otimes
\left(
\begin{array}{cc}
1 & 0\\
0 & 0
\end{array}
\right)
+ i
\left(
\begin{array}{cc}
0 & 0\\
0 & 1
\end{array}
\right)
\otimes
\left(
\begin{array}{cc}
0 & 0\\
0 & 1
\end{array}
\right).
\end{equation}
Then $A$ is normal matrix with eigenvalues $0,1,i$. It is easy to see that 
$1 \in \ProductNumRange{A}$ and $i \in \ProductNumRange{A}$, but 
$(1+i)/2 \not\in\ProductNumRange{A}$.
Actually, by direct computation we have
\begin{equation}
\ProductNumRange{A}= \left\{ x+yi : 0\leq x, 0\leq y, \sqrt{x}+\sqrt{y}\leq 1\right\}.
\end{equation}
Product numerical range of matrix $A$ is presented in Figure~\ref{fig:example-nonconvex}.
\end{example}
\halmos

Product numerical range forms a~nonempty set for a~general operator. 
In particular it contains the barycenter of the spectrum.

\begin{property}\label{barycentergeneral}
Product numerical range of $A\in \MM_{KM}$ 
includes the barycenter of the spectrum, 
\begin{equation}
\frac{1}{KM}\; {\tr} A \ \in \  \ProductNumRange{A}.
\end{equation}
\end{property}
\proof
Let $A$ be an operator acting on a~tensor product Hilbert space
$\mathcal{H}_{K}\otimes\mathcal{H}_{M}$. Let us write
\begin{equation}
\label{barycenterreduct}
 \frac{1}{KM}\; \tr A=
\frac{1}{KM}\; \tr\left(A(\1\otimes\1)\right)=\frac{1}{K}\sum_{i=1}^M\frac{1}{M}
\tr\left(A\left(\1\otimes\ketbra{\psi_i}{\psi_i}\right)\right),
\end{equation}
where $\left\{\psi_i\right\}_{i=1}^M$ is an arbitrary orthonormal basis in $\mathcal{H}_M$. 
The last sum in \eqref{barycenterreduct} is a~convex combination of elements in the numerical 
range of $\tr_1 A$, where $\tr_1$ denotes the partial trace with respect to $\mathcal{H}_K$. 
Remember that $\LambdaRm\left(\tr_1 A\right)$ is convex. Hence there exists an element 
$\psi\in\mathcal{H}_M$ of norm one such that 
\begin{equation}
\sum_{i=1}^M\frac{1}{M}
\tr\left(A\left(\1\otimes\ketbra{\psi_i}{\psi_i}\right)\right)
=\bra{\psi}\tr_1 A \ket{\psi}=\tr\bigl(A\left(\1\otimes\ketbra{\psi}{\psi}\right)\bigr).
\end{equation}
By repeating the same trick, we can replace the remaining identity with a~single
one-dimensional projector and obtain
\begin{eqnarray}\label{elimidentities}
\frac{1}{KM}\;  \tr A 
&=& \frac{1}{K} \tr \left(A\left(\1\otimes\ketbra{\psi}{\psi}\right)\right) \\
\nonumber
&=& \tr\left(A\left(\ketbra{\phi}{\phi}\otimes\ketbra{\psi}{\psi}\right)\right)
=\bra{\phi\otimes\psi}A\ket{\phi\otimes\psi},
\end{eqnarray}
for some $\ket{\phi}\in\mathcal{H}_K$ and $\ket{\psi}\in\mathcal{H}_M$ of norm
one. The last equality in \eqref{elimidentities} means the same as 
$ (\tr A)/KM \; \in \;  \ProductNumRange{A}$, which we wanted to prove.
\halmos

In the particular case $\tr A=0$, Property \ref{barycentergeneral} has already been used
in~\cite{DFY08}. The above reasoning can be generalized to the multipartite
case (cf. Section \ref{sec:multipartite}).

Note that the barycenter does not have to lie in the interior of the product numerical range. For example, point $\left(\frac{1}{4},\frac{1}{4}\right)$, denoted by the black cross in Figure \ref{fig:example-nonconvex}.

\begin{property}
Product numerical radius is a~vector norm on matrices, but it is not a~matrix
norm. Product numerical radius is invariant with respect to local unitaries,
which have the tensor product structure.
\end{property}
\proof
The only hard part of the proof is positivity. In the Hermitian case this
follows from Proposition~\ref{prop:existence}, the non-Hermitian case can be
reduced to Hermitian by the projection property.
\halmos

\begin{property}\label{prop:diagonal-parametrization2x2}
If $A\in \MM_{KM}$ can be diagonalized to $\Sigma$ using product unitary 
matrices $U\in \MM_K, V\in \MM_M$ (\ie{} there exist unitary $U, V$ such that 
$\left(U\otimes V\right) A \left(U^\dagger \otimes V^\dagger\right)=\Sigma$) then 
\begin{equation}
\ProductNumRange{A}=\{z:z=\left((p_1,p_2,\ldots, p_K)\otimes (q_1,q_2,\ldots, q_M)\right)\cdot
(\Sigma_{1,1},\Sigma_{2,2},\ldots,\Sigma_{KM,KM})\},
\end{equation}
where $\sum_k p_k=1$, $\sum_m q_m=1$ and $p_k, q_m\geq 0$.
\end{property}
\proof This follows from more general 
Proposition~\ref{prop:diagonal-parametrization} in Section~\ref{sec:multipartite}.
\halmos

\begin{example} 
We can apply the last property to Example \ref{example-nonconvex}.
The matrix in this example is given by
\begin{equation} 
A = \diag(1,0,0,i).
\end{equation} 
By Property \ref{prop:diagonal-parametrization2x2} we have 
the following parametrization of the product numerical range of $A$:
\begin{equation}
p q + i (1 - p) (1 - q), \; \; p,q \in [0,1].
\end{equation}
\end{example}
\subsubsection{Relation to Minkowski geometric algebra}
\label{sec:tens_prod}
We shall start this section by recalling the Minkowski geometric algebra of complex sets 
as developed by Farouki \etal \ \cite{FMR01}.
For any sets $Z_1$ and $Z_2$  in the complex plane
one defines their \emph{Minkowski sum},

\begin{equation}
Z_1  \mplus Z_2 =  
\left\{z: z=z_1+z_2 , \ z_1 \in Z_1,  z_2 \in Z_2 \right\} 
\label{Msum} 
\end{equation}
and  \emph{Minkowski product},
\begin{equation}
Z_1 \mprod  Z_2 =  
\left\{z: z=z_1 z_2, \ z_1 \in Z_1, \ z_2 \in Z_2 \right\} .
\label{Mprod} 
\end{equation}
Note that the above operations are not denoted by $\oplus$ and $\otimes$
as in the original paper \cite{FMR01},
in order to avoid the risk of confusion with the direct sum of operators or their
tensor product.

A simple lemma concerning the Minkowski sum and product has interestingly deep
consequences. Let us define the Kronecker sum of two operators as 
\begin{equation}
A \oplus B =  A \kron \1  +  \1 \kron B.
\end{equation}

\begin{lemma} \label{lemma:KronToMin}
Product numerical range of the Kronecker product of 
arbitrary operators is equal to the Minkowski product 
of the numerical ranges of both factors,
\begin{equation}
  \label{rangetensorpr}
\ProductNumRange{A \otimes B}   =  \NumRange{A} \mprod  \NumRange{B},
\end{equation}
while product numerical range of the Kronecker sum of 
arbitrary operators is equal to the Minkowski sum
of the numerical ranges of both factors,
\begin{equation}\label{rangetensorsum}
\ProductNumRange{A \oplus B}    =  \NumRange{A} \mplus  \NumRange{B}.
\end{equation}
\end{lemma}

Proof follows directly from the definition of the product numerical range. First
part of the above lemma has already been used in \cite{dirr08relative}. Observe
that the definition (\ref{Mprod}) and the property (\ref{rangetensorpr}) can be
naturally generalized to an arbitrary number of factors.

Thus the problem of finding the product numerical range of a~tensor product 
can be analyzed by checking what sorts of subsets of the complex plane 
one can obtain by multiplying two or more numerical ranges.
This very problem has recently been investigated in a~series of papers
by Farouki \etal\ -- see \cite{FMR01,FP02} and references therein.
For instance, the structure of the Minkowski product of several discs
in the complex plane was analyzed in detail by Farouki and Pottmann~\cite{FP02}.

The above results concerning Minkowski product can be used directly to find
the product numerical range of a~tensor product of an arbitrary number of
factors acting on two-dimensional subspaces. The numerical range of any matrix
of order two forms an ellipse \cite{Do57}, which may degenerate to an interval
or a~point. For instance it is known that the numerical range of the matrix
$  X=\left( 
\begin{smallmatrix}
    c &\ 2r \\
    0 & c 
\end{smallmatrix} \right)$,
forms a~disk of radius $|r|$ centered at $c$. 
Consider a~family of operators with a~tensor product structure 
\begin{equation}
\label{proddisks}
Y(r_1,r_2) = 
X_1 \otimes X_2 =
\left(
\begin{array}{cc}
1 & \ 2 r_1\\
0 & 1
\end{array}
\right)
\otimes
\left(
\begin{array}{cc}
1 & \  2r_2\\
0 & 1
\end{array}
\right).
\end{equation}
\begin{figure}[h]
 \includegraphics[width = \textwidth ]{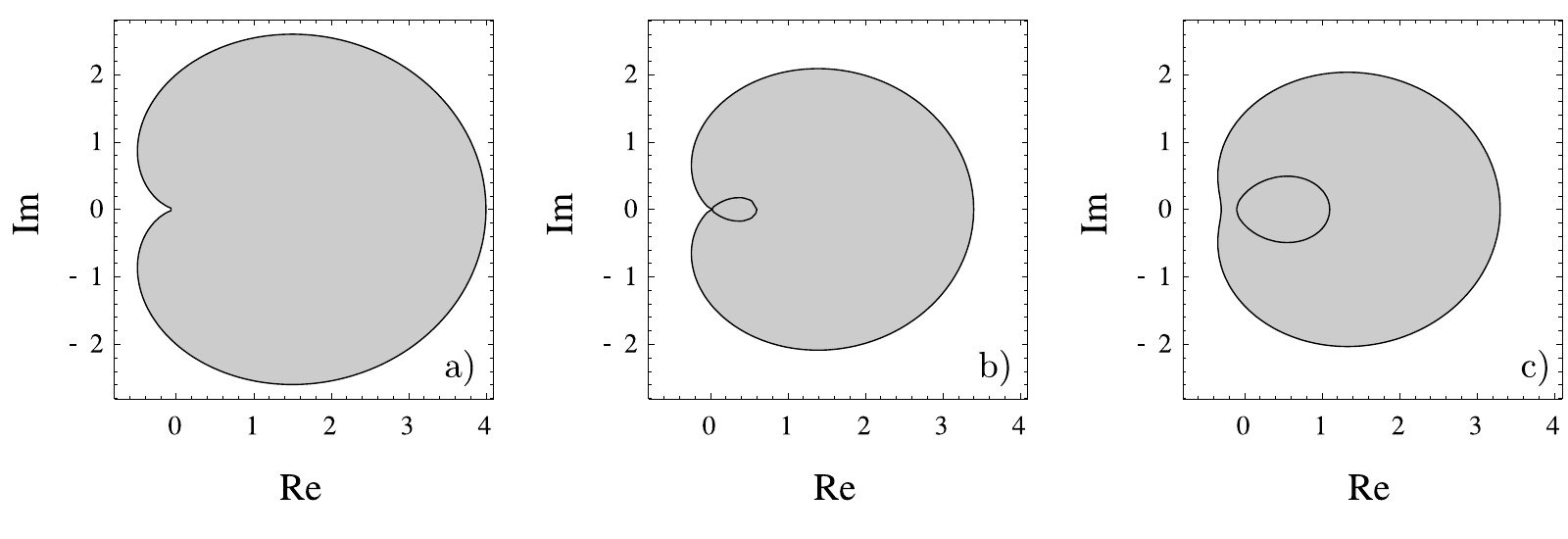}
 \caption{Product numerical range for the operator $Y$ 
defined in (\ref{proddisks})
with $(r_1,r_2)$ equal to a) $(1,1)$ (cardioid), 
  b) $(0.7, 1)$ (lima{\c c}on of Pascal), and  
  c) $(0.5,1.2)$ (Cartesian oval).}\label{fig:cardioid}
\end{figure}

The product numerical range of $Y$ takes different shapes depending 
on the values of the radii $r_1$ and $r_2$ of both discs.
According to the results of \cite{FP02}, one may find the values of these parameters 
for which the boundary of the product numerical range is a~cardioid,
a~lima\c{c}on of Pascal, or the outer loop of a~Cartesian oval -- 
see Figure~\ref{fig:cardioid}.

Analysis of the Minkowski product of two sets becomes easier if none of them
contains $0$. In such a~case one may use the log-polar coordinates in the
complex plane and reduce the Minkowski product to a~Minkowski sum in the new
coordinates.

Let us now consider the opposite case.
\begin{lemma}
\label{lemmastar}
If the numerical range of one of the factors $A_1$ or $A_2$ contains $0$,
then the product numerical range of the tensor product 
 $\ProductNumRange{ A_1\otimes A_2 }$ is star-shaped.
\end{lemma}
\proof
We have
\begin{equation}
\label{starsum}
\ProductNumRange{A_1\otimes A_2} =  
\bigcup_{z\in\NumRange{A_2}}z \NumRange{A_1}.
\end{equation}
Without loss of generality, we may assume that $0 \in \NumRange{A_1}$. 
It is known that the numerical range $\LambdaRm (A_1) $ is convex.
Hence $z\NumRange{A_1}$ is star-shaped with respect to $0$ for 
arbitrary $z \in \Cplx$. 
We get that $\ProductNumRange{A_1\otimes A_2}$ is star-shaped with respect 
to $0$ and therefore simply connected.
\halmos

It is conceivable that the assertion of Lemma \ref{lemmastar} may also hold
without the assumption $0\in\Lambda\left(A_1\right)\cup\Lambda\left(A_2\right)$,
but so far, we were not able to prove or disprove this.

We are now in a~position to formulate the main result of this section.
\begin{proposition}
\label{productoftwo}
Let $A_1,A_2$ be arbitrary operators. The product numerical range of $A_1\otimes
A_2$ is simply connected.
\end{proposition}
The above proposition, proved in \ref{appendix_b}, does not hold for a
three-fold tensor product. In this way we confirm the conjecture formulated in
\cite{dirr08relative} that one needs to work with at least tripartite systems to
construct a~tensor product operator whose product numerical range is not simply
connected.

\subsubsection{Inclusion properties}
\label{sec:bipartite}

Generically, operators acting on a~bipartite Hilbert space~(\ref{HKM}) do not
exhibit the tensor product form. However, treating operators as vectors in the
Hilbert-Schmidt space endowed with the Hilbert-Schmidt scalar product,
$\scalar{A}{B} = {\rm Tr}A^{\dagger} B$, we may use the operator Schmidt
decomposition. In close analogy to regular Schmidt decomposition, any operator $X$
acting on ${\cal H}_K \otimes {\cal H}_K$ can be decomposed as a~sum of not more
than $K^2$ terms,
\begin{equation}
\label{opschmidt}
X = \sqrt{\mu_1} A_1\kron B_1 + \dots + \sqrt{\mu_{K^2}} A_{K^2}\kron B_{K^2} .
\end{equation}

To find the explicit form of this decomposition, it is convenient to use the
reshuffled matrix, $Y=X^R$, such that in the product basis it
consists of the same entries as the original matrix $X$, but ordered
differently, $\bra{i,j}Y\ket{ i' , j'} =\bra{i,i'}X\ket{j , j'}$. 
Then the non-negative components $\mu_i$ of the Schmidt vector are
equal to the singular values of the non-negative matrix $YY^{\dagger}$ of order
$K^2$, while operators $A_i$ and $B_i$ with $i=1,\dots K^2$ are obtained by
reshaping eigenvectors of Hermitian matrices, $YY^{\dagger}$ and $Y^{\dagger}Y$,
respectively \cite{BZ06}.

Making use of the Schmidt decomposition (\ref{opschmidt}) of an arbitrary
bipartite operator $X$ we can formulate a~proposition concerning its product
numerical range.
\begin{proposition}
\begin{eqnarray}
\ProductNumRange{X}  &\subset&  
 \sqrt{\mu_1} \left(\NumRange{A_1} \mprod  \NumRange{B_1} \right)
  \mplus  \dots \mplus  
\sqrt{\mu_{K^2}}  \left(\NumRange{A_{K^2}} \mprod \NumRange{B_{K^2}}\right) 
\nonumber \\
&&=\ProductNumRange{\sqrt{\mu_1}  A_1 \otimes B_1}  \mplus \dots \mplus  
\ProductNumRange{\sqrt{\mu_{K^2}} A_{K^2} \otimes B_{K^2}} .
\label{lnrschmidt}
\end{eqnarray}
\end{proposition}
\proof 
It follows directly from Property \ref{prop:subadditivity} and from the fact
that the product numerical range of a~Kronecker product of operators equals the
Minkowski product of their numerical ranges (Lemma \ref{lemma:KronToMin}).  
\halmos

\subsection{Hermitian case}
\label{sec:lr-Hermitian}

\subsubsection{Basic properties}
Consider a~Hermitian operator $X=X^{\dagger}$ with ordered spectrum $\lambda_1
\le \lambda_2 \le \dots \le \lambda_N$. A~closed interval
$[\lambda_k,\lambda_{k+1}]$ will be called a~\emph{segment} of the spectrum.
Thus the spectrum consists of $N-1$ segments, some of which reduce to a~single
point, if the spectrum is degenerated. The numerical range of a~Hermitian
operator reads $\NumRange{X}=[\lambda_1,\lambda_N]$~\cite{HJ2}. In this case
also product numerical range is given by an interval,
$\ProductNumRange{X}=[\lambdaProdMin,\lambdaProdMax]$, bounded by the points
$\lambdaProdMin$ and $\lambdaProdMax$ -- the extremal points of product
numerical range under the set of all product states. Hence, these points
determine the maximal (the minimal) expectation values of an observable $X$
among all product pure states.

Observe that these extremal values determine the product numerical radius,
\begin{equation}
 r^{\otimes}(X) ={\rm max}\{ |\lambdaProdMin(X)|, |\lambdaProdMax(X)| \}.
\end{equation}
This relation simplifies for quantum states. Their positivity implies that 
$r^{\otimes}(\rho) = \lambdaProdMax(\rho)$.

Note that the product numerical range by definition depends on the concrete form
of the tensor product Hilbert space $\mathcal{H}_N=\mathcal{H}_K \otimes
\mathcal{H}_M$. For instance, if $X$ also possesses the similar product
structure, $X=X_A \otimes X_B$, then all its eigenstates, $\ket{\phi_i},
i=1,\dots,N$, are of the product form and are called \emph{separable states}.
Any pure state which is not of the product form is called \emph{entangled\/}.
Thus in this particular case both ranges are equal,
$\ProductNumRange{X}=\NumRange{X}=[\lambda_1,\lambda_N]$. This is also the case
if the eigenstates $\ket{\phi_1}$ and $\ket{\phi_N}$, corresponding to the
extremal eigenvalues $\lambda_1$ and $\lambda_N$, are of the product form. Hence
the product structure of $X$ is not necessary to assure that both ranges
coincide.

Let us consider an arbitrary orthonormal product basis,
 $\ket{i,j} =\ket{i} \otimes \ket{j}$ of $\mathcal{H}_N$.
The states $ \ket{i} \in {\cal H}_K$, $i=1,\dots, K $ and
 $\ket{j} \in {\cal H}_M$, $i=1,\dots, M$ 
satisfy the orthogonality relation
 $\scalar{i,j}{i',j'}=\delta_{i,i'} \delta_{j,j'}$.
Let us also use the composed indices denoted by Greek letters, 
$\mu=(i,j), \nu=(i',j')$,
to represent $X$ in the product basis
\begin{equation}
 X_{\mu\nu} =   \bra{ \mu }X\ket{\nu} =
 \bra{i,j}X\ket{i' , j'}.
 \label{xmn}
\end{equation}
In such a~product representation any diagonal element $X_{\mu\mu}$ denotes the 
expectation value of $X$ in a~product state $\ket{\mu}=\ket{i,j}$, so
it belongs to the product numerical range, $X_{\mu\mu} \in \ProductNumRange{X}$.

Interestingly, some information about the product numerical range can be
obtained even without specifying a~concrete tensor product structure in
$\mathcal{H}_N$ and a~product basis.

\subsubsection{Invariant features}
In this section we investigate some basic properties of the product numerical
range of a~Hermitian operator, $X = X^{\dagger}$, which hold independently from
the splitting of the Hilbert space $\mathcal{H}_N=\mathcal{H}_K \otimes
\mathcal{H}_M$. In the general case of an operator $X$ acting on a~$K\times M$
space, not so much can be said about its product range.

\begin{proposition}\label{prop:existence} 
Product numerical range forms, by definition, a~non-empty subset of the
numerical range, $\emptyset \neq \ProductNumRange{X} \subset \NumRange{X}$.
Furthermore, $\ProductNumRange{X}$ reduces to a~single point $\lambda$ if and
only if the operator $X$ is proportional to the identity, $X=\lambda \1$.
Speaking in terms of the standard Lebesgue volume of an interval we arrive at
the following statement. If {\rm Vol}$[\NumRange{X}] >0$, then {\rm
Vol}$[\ProductNumRange{X}]>0$.
\end{proposition}

The assertion that $\ProductNumRange{X}$ is a single point if and only if $X=\lambda\1$ holds also when we do not assume Hermiticity of $X$.

\begin{proposition}\label{prop:convex-hermit} 
For any Hermitian operator $X$, acting on an $N$-dimensional Hilbert space, its
product numerical range $\ProductNumRange{X}$ is convex and forms an interval of
the real line.
\end{proposition}

\proof{
Let us assume that $a_1$ and $a_2$ belong to $\ProductNumRange{X}$. Hence there exist 
two pairs of product vectors, such that $a_1=\bra{x_1, y_1} X \ket{x_1,y_1}$ 
and $a_2=\bra{x_2, y_2} X \ket{x_2,y_2}$. Since it is possible to find
a~continuous family of vectors $\ket{\phi_A (s)} \in {\cal H}_K$, which 
interpolates between $\ket{x_1}$ and  $\ket{x_2}$ and another family 
$\ket{\phi_B(s)} \in \mathcal{H}_M$, which interpolates between $\ket{y_1}$ 
and $\ket{y_2}$, the expectation values of $X$ between these product states 
interpolate between $a_1$ and $a_2$. Thus the entire interval belongs to the 
product numerical range, $[a_1,a_2] \subset \ProductNumRange{X}$.}
\halmos

Note that the above reasoning, applied to an arbitrary (non-Hermitian) operator
$X$, shows that the product numerical range forms a~connected set in the complex
plane. However, as shown in further sections of this paper, this set does not
need to be convex nor simply connected.

Before we turn to more specific theorems, let us invoke a~useful result. 

\begin{lemma}
\label{lem:min-subspace} 
Consider a~tensor product complex Hilbert space $\mathcal{H}_N=\mathcal{H}_K \otimes 
\mathcal{H}_M$. Then we have the following:
\begin{enumerate}
 \item[\emph{a)}] any subspace $S_d\subset \mathcal{H}_N$ of dimension 
$d=(K-1)(M-1)+1$ contains at least one separable state,\label{min-subspace-part1}
 \item[\emph{b)}] there exists a~subspace of dimension $d-1=(K-1)(M-1)$ which
does not contain any separable state.
\end{enumerate}
\end{lemma}

\proof
Part a) of the above Lemma follows directly from Proposition 6 in a~paper 
by Cubitt \etal\ \cite{CMW08}, but this fact was earlier proven 
in~\cite{Wa02,Pa04}.

Concerning part b), 
for given integers $K,M$ we give an example of 
family of $(K-1)\times (M-1)$ matrices $A^{(ij)}$
$i=1,\ldots,K-1$ and $j=1,\ldots,M-1$ such that no linear combination of 
$\left\{A^{(ij)}:1\leq i < K, 1\leq j < M\right\}$ is of rank one.

 Let $A^{(ij)} = E^{(ij)}+E^{(i+1,j+1)}$ for $1\leq
i < K$ and $1\leq j < M$. Explicitly, let us denote by $E^{(ij)}$ the $K\times M$ matrix containing 1 at the position
$(i,j)$ and zeros elsewhere. Suppose that $T=\sum_{i=1}^{K-1} \sum_{j=1}^{M-1}\alpha_{ij} 
A^{(ij)}$ is of rank one. As the top-right entry of $T$ is zero, then the
first row or the $M$-th column of $T$ is zero. Thus we have $\alpha_{1,j}=0$ for
all $j$ or $\alpha_{i,M}=0$ for all $i$. Deleting the first row or the $M$-th
column of $T$ we can proceed by induction on $K$ and $M$ to get the desired
result.

Thus the subspace spanned by the vectorizations 
$\sum_{k,l}A^{\left(ij\right)}_{lk}\left|k,l\right>$ of the matrices $A^{(ij)}$ 
does not contain any product state. 
\halmos

A~subspace containing  no product states is called \emph{completely entangled}
and it was shown in \cite{WS08} that a~generic subspace of dimension $(K-1)(M-1)$
possesses this property.

In the proof of part b) of Lemma \ref{lem:min-subspace}, we used the fact that any
pure state in a~$N=KM$ dimensional bipartite Hilbert space can be represented in terms
its \emph{Schmidt decomposition},
\begin{equation}
   \ket{\psi}   =   \sum_{i=1}^K \sum_{j=1}^M 
   A_{ij} \ket{i} \otimes \ket{j}
      =  
  \sum_{i=1}^{K} \sqrt{\mu_i}\; \ket{i'} \otimes \ket{i''}.
  \label{Aij}
  \end{equation}
We have assumed here that $K\le M$ and denoted a~suitably rotated product basis
by $\ket{i'} \otimes \ket{i''}$. The eigenvalues $\mu_i$ of a~positive
matrix $AA^{\dagger}$ are called the \emph{Schmidt coefficients} of the
bipartite state $\ket{\psi}$. The normalization condition 
$|\psi|^2=\scalar{\psi}{\psi}=1$ implies that 
$||A||^2_{\rm HS}=\tr AA^{\dagger}=1$, so the Schmidt coefficients $\mu_i$ form 
a~probability vector.

The state $\ket{\psi}$ is separable iff the $K \times M$ matrix of coefficients
$A$ is of rank one, so the corresponding vector of the Schmidt coefficients is
pure. Hence Lemma \ref{lem:min-subspace} is equivalent to the following. 

\begin{lemma}
Consider a~set of $k$ complex rectangular matrices $A_i$ of size $K \times M$, 
which are orthogonal with respect to the Hilbert-Schmidt scalar product, 
$\scalar{A_i}{A_j} = {\rm Tr}A_i^{\dagger}A_j=\delta_{ij}$. 
\renewcommand{\labelenumi}{\alph{enumi})}
\begin{enumerate}
\item 
 If $k \ge d=(K-1)(M-1)+1$ there exists a~complex vector $\vec c$ determining
a~linear combination of these matrices, $A_{\rm av}=\sum_{i=1}^k c_i A_i$, such
that $A_{\rm av}$ is of rank one.
\item 
 Moreover, simple dimension counting arguments presented in~\cite{CMW08} imply
that, for a~generic set of $(d-1)$ such matrices $A_i$, a~rank one linear
combination does not exist.
\end{enumerate}
\end{lemma}

\begin{proposition}[General $K \times M$ problem] 
 \label{prop:KtimesM}  
The following lower bound for the product numerical radius is true,
 $$\lambdaProdMax \ge \lambda_{K+M-1}$$
 and a~symmetric upper bound for the product minimum holds, 
 $$\lambdaProdMin \le \lambda_{(K-1)(M-1)+1}.$$

Furthermore,
there exist $X_1=X_1^{\dagger}$ and $X_2=X_2^{\dagger}$ 
acting on ${\cal H}_{KM}$ such that  
$\lambdaProdMax(X_1) < \lambda_{K+M}(X_1)$ and 
$\lambdaProdMin(X_2) > \lambda_{(K-1)(M-1)}(X_2)$.
\end{proposition}

Before proving this proposition let us state some of its implications.
For any Hermitian operator $X$ acting on the  $2 \times 2$ Hilbert space
its product numerical range contains the central segment of the spectrum,
$[\lambda_2,\lambda_3] \subset \Lambda^{\otimes}(X)$.
A similar statement holds for the $2 \times m$ problem, 
$[\lambda_m,\lambda_{m+1}] \subset \Lambda^{\otimes}(X)$.
Proposition \ref{prop:KtimesM} implies that 
for any $X$ acting on a~$3 \times 3$ system
 $\lambda_5 \in \Lambda^{\otimes}(X)$.

\proof{
Lemma \ref{lem:min-subspace} implies in this case that any $d=(M-1)(K-1)+1$
dimensional subspace contains a~product state. Applying this to the subspace
spanned by the $(M-1)(K-1)+1$ eigenstates corresponding to the $(M-1)(K-1)+1$
smallest eigenvalues we obtain a~bound $\lambda_{\min}^{\otimes} \le
\lambda_{(M-1)(K-1)+1}$. The $(M-1)(K-1)+1$-dimensional subspace spanned by the
eigenstates corresponding to largest eigenvalues also contains at least one
product state, which implies that $\lambdaProdMax \ge \lambda_{K+M-1}$.

To demonstrate the optimality, it is enough to find a~single
operator $X$ of size $KM$ such that $\lambda_{(K-1)(M-1)} \notin
\ProductNumRange{X}$. By Lemma \ref{lem:min-subspace}, there exists a~subspace
$\mathcal{K}$ of dimension $(K-1)(M-1)$ which does not contain a~separable
state. Let $\ket{\psi_1},\dots,\ket{\psi_{(M-1)(K-1)}}$ be some orthonormal
basis of this subspace, and let vectors
$\ket{\psi_{(M-1)(K-1)+1}},\dots,\ket{\psi_{KM}}$ enlarge the previous system to
the orthonormal basis of $\mathcal{H}_{KM}$. We define matrix $X$ as
\begin{equation}
 X = \sum_{i=(M-1)(K-1)+1}^{KM} \ketbra{\psi_i}{\psi_i}.
\end{equation}
Each product unit length vector can be written in the above basis as
\begin{equation}
\ket{x,y} =  \sum_{i=1}^{(M-1)(K-1)} \alpha_i \ket{\psi_i} + \sum_{i=(M-1)(K-1)+1}^{KM} \beta_i \ket{\psi_i}.
\end{equation}
Because the subspace $\mathcal{K}$ does not contain any separable state, in the
above decomposition of a~product state there exists $i \in
\{(M+1)(K-1)+1,\dots,KM\}$ such that $\beta_i \neq 0$. Thus for any product
state we have
\begin{equation}
\bra{x,y} X \ket{x,y} \  > \  0 \ = \  \lambda_{(M-1)(K-1)}.
\end{equation}
}\halmos

Note that it is not possible to obtain complete
information about the position of the product numerical range relying only on
unitarily invariant properties. In general, one has to consider a concrete splitting of the composite Hilbert space.

\section{Generalizations}
\label{sec:generalization}

\subsection{Multipartite operators}
\label{sec:multipartite}

\subsubsection{Definition}
It is straightforward to generalize the notions of product numerical range 
to spaces with a~multipartite structure. In this
section we study Hilbert spaces formed by a~tensor product of an arbitrary
number of factors, 
\begin{equation}
\mathcal{H}_N \ = \  \mathcal{H}_{m_1} \otimes \cdots \otimes \mathcal{H}_{m_k},
\end{equation}
with $N=m_1\dots m_k$. 

For instance, in the definition of product numerical radius 
one has to perform the maximization over the direct product group of
\emph{local unitary transformations},
\begin{equation}
\label{localuni}
U_{\rm loc}\ := \ U(m_1) \times U(m_2) \times \cdots U(m_k) \ ,  
\end{equation}
embedded in the group of global unitaries, $U(N)=U(m_1\dots m_k)$.

\subsubsection{Parametrization}

The following proposition gives a~parametrization of the product numerical range
for operators which can be diagonalized by a product of unitary transformations. 
\begin{definition}[Product diagonalizable operator]
We call operator $A\in \MM_{m_1 m_2 \ldots m_k}$ 
product diagonalizable iff there exists unitary operators 
$U_{m_1}\in \MM_{m_1}, U_{m_2}\in \MM_{m_2}, \ldots, U_{m_k}\in \MM_{m_k}$ such 
that 
\begin{equation}
(U_{m_1}\otimes U_{m_1} \otimes \ldots \otimes U_{m_k}) A (U_{m_1}\otimes U_{m_1} \otimes \ldots \otimes U_{m_k})^\dagger=\Sigma
\end{equation}
and $\Sigma$ is a~diagonal matrix.
\end{definition}
In the special case where the operator $A$  is Hermitian and positive,
such states are called {\it classically correlated} 
\cite{OHHH02} or {\it locally diagonalizable} \cite{BZ06}.

\begin{proposition}\label{prop:diagonal-parametrization}
Let $A$ be a~product diagonalizable operator on  
$\mathcal{H}_{m_1 m_2 \dots m_k } =  \mathcal{H}_{m_1} \kron \mathcal{H}_{m_2} \kron \ldots \kron \mathcal{H}_{m_k}$.
By 
\begin{equation} 
\lambda_{l_1,l_2,\dots, l_k} = \bra{l_1,l_2,\dots , l_k} A \ket{l_1,l_2,\dots  , l_k}
\end{equation} 
we denote the vector of diagonal elements of $A$ after diagonalization using
product unitary matrices. With the above notation, we have the parametrization
of the product numerical range of $A$,
\begin{equation} \label{eqn:parametrizationOfLNR}
\ProductNumRange{A} = 
\left\{z : z =  
 \sum_{l_1=0}^{m_1-1} \sum_{l_2=0}^{m_2-1} \dots \sum_{l_k=0}^{m_k-1} 
p^{(1)}_{l_1} p^{(2)}_{l_2} \dots p^{(k)}_{l_k} \lambda_{l_1,l_2,\dots l_k} \right\},
\end{equation}
 where $p^{(r)}_0,\dots,p^{(r)}_{m_r-1} \geq 0$ and 
$ p^{(r)}_0 + \dots + p^{(r)}_{m_r-1}= 1 $
for $r=1, \dots, k$.
\end{proposition}
\begin{remark}
Note that if $A$ is a~product diagonalizable operator on  
$\mathcal{H}_{2^k} =  \mathcal{H}_{2} \kron \mathcal{H}_{2} \kron \ldots \kron \mathcal{H}_{2}$,
then the above parametrization simplifies and we have
\begin{eqnarray} \label{eqn:parametrization2^kOfLNR}
\nonumber
\ProductNumRange{A} = 
\Big\{z : z &=&  
\left( \left\{ p^{(1)},1-p^{(1)} \right\} \kron  \dots \kron \left\{p^{(k)},1-p^{(k)}\right\} \right)
\cdot \lambda ,\\
&&
p^{(1)},\dots,p^{(k)} \in [0,1]  
\Big\},
\end{eqnarray}
where '$\cdot$' denotes the scalar product.
\end{remark}

A~proof of the above proposition is given in~\ref{sec:proof-th-diagonal-parametrization}. 
Later we will use this proposition to study some concrete examples.

\subsubsection{Properties}
\begin{proposition}\label{barycentermultigeneral}
Product numerical range of an operator $X$ acting on a~tensor product Hilbert space 
$\mathcal{H}_N = \mathcal{H}_{m_1} \otimes \cdots \otimes \mathcal{H}_{m_k}$ 
includes the barycenter of the spectrum \ie, 
\begin{equation}
\frac{1}{N}{\tr} X\; \in \; \ProductNumRange{X}\; .
\end{equation}
\end{proposition}
\proof
We can follow the same line of reasoning as in the proof of Property \ref{barycentergeneral}.
\halmos

The next example, adopted from \cite{thomas08significance}, shows that, in the
case of a~three--partite system, product numerical range is not necessarily
simply~connected.

\begin{figure}[ht!]
\begin{center}
\label{fig:nonsimply-example}
\subfigure[$\LambdaRm^{\otimes}$ of matrix given by~(\ref{eqn:thomas-nonsimply-example}).]{
	\label{fig:non-simply-connected}
	\includegraphics[width=0.45\textwidth,scale=0.6]{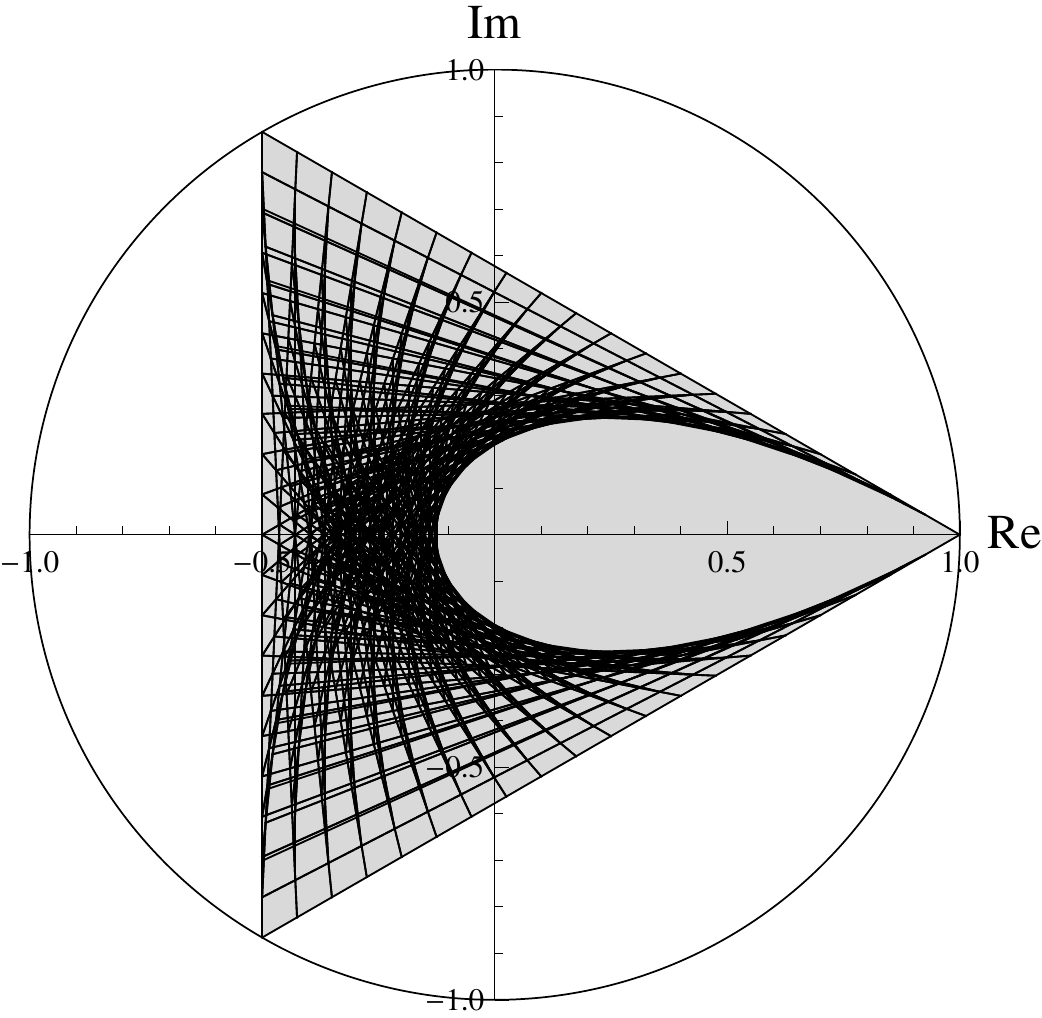}
}
\subfigure[$\LambdaRm^{\otimes}$ of matrix given by~(\ref{eqn:thomas-simply-example}).]{
	\label{fig:simply-connected}
	\includegraphics[width=0.45\textwidth,scale=0.6]{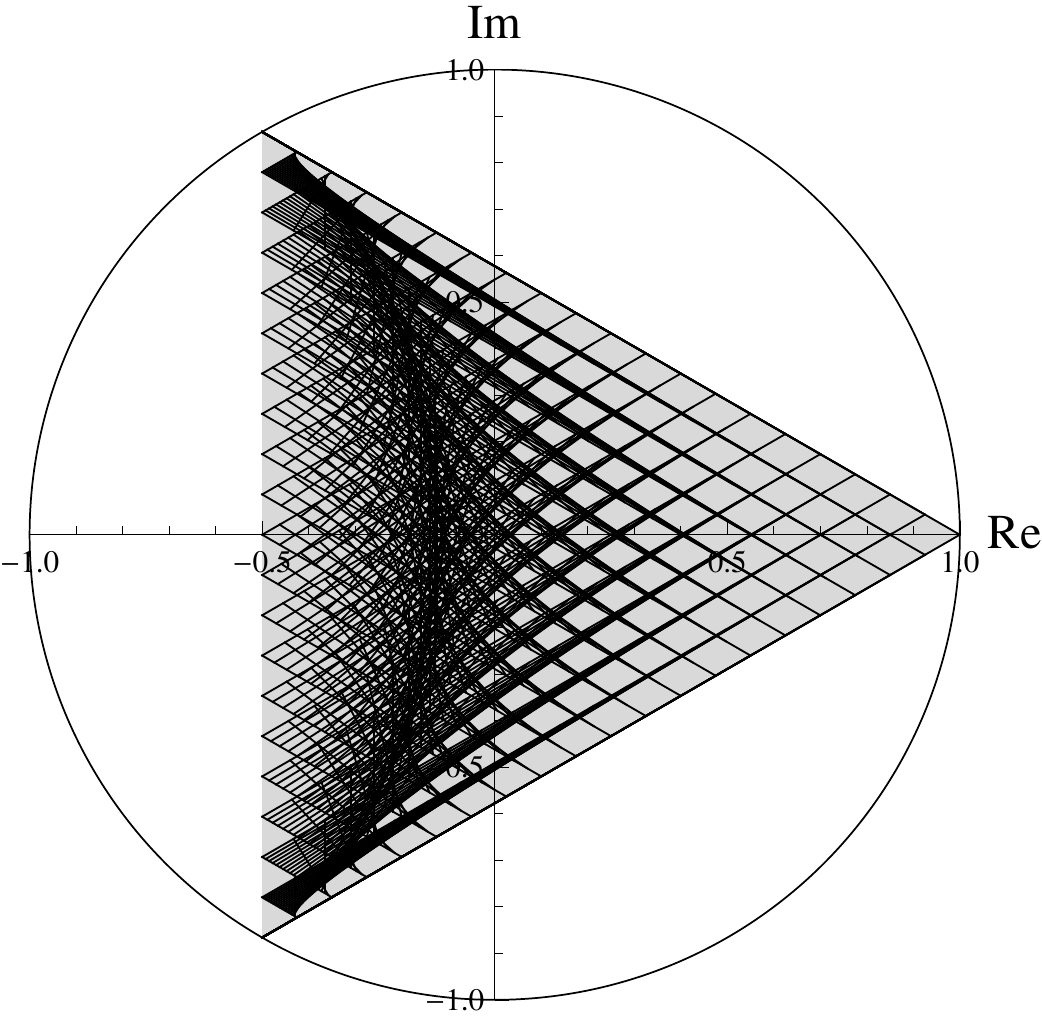}
}
\caption{Figures obtained from parametrizations (\ref{eqn:par-U1}) and (\ref{eqn:par-U2}). 
In both cases the numerical range 
 forms an equilateral triangle plotted in grey.}
\end{center}
\end{figure}

\begin{example}
Consider two unitary matrices $U_1$ and $U_2$ 
written in the standard computational basis
$\{\ket{000}, \ket{001}, \dots , \ket{111} \}$,
\begin{equation}
\label{eqn:thomas-nonsimply-example}
U_1 =\diag(1,e^{\frac{2 i \pi }{3}},e^{\frac{2 i \pi }{3}},e^{-\frac{2 i \pi }{3}},
e^{\frac{2 i \pi }{3}},e^{-\frac{2 i \pi }{3}},e^{-\frac{2 i \pi }{3}},1)
\end{equation}
and
\begin{equation}
\label{eqn:thomas-simply-example}
U_2 =\diag(1,e^{\frac{2 i \pi }{3}},e^{\frac{2 i \pi }{3}},e^{\frac{2 i \pi }{3}},
e^{-\frac{2 i \pi }{3}},e^{-\frac{2 i \pi }{3}},e^{-\frac{2 i \pi }{3}},1).
\end{equation}
\rm Both $U_1$ and $U_2$ have identical eigenvalues, however they have different
eigenvectors associated with these values.

Product numerical range of the operator $U_1$ acting on a~three--partite system
is not simply connected, as shown in Fig.~\ref{fig:non-simply-connected}. 
However, exchanging the position of two middle eigenvalues one obtains
a unitary operator $U_2$ with the same spectrum, for which the product numerical
range is convex and coincides with the standard numerical range,
 $\ProductNumRange{U_2}=\NumRange{U_2}$, 
see: Fig.~\ref{fig:simply-connected}.


Making use of Proposition \ref{prop:diagonal-parametrization}, we can 
explicitly parameterize the 
product numerical range of these matrices,
\begin{equation} \label{eqn:par-U1}
\ProductNumRange{U_1} =
\Big\{ z : z = 
\left(
\{p,1-p\} \kron \{q,1-q\} \kron \{r,1-r\}
\right) \cdot 
(\{U_1\}_{ii})_{i=1}^8 \Big \}
\end{equation}
and
\begin{equation}\label{eqn:par-U2}
\ProductNumRange{U_2} =
\Big\{ z : z = 
\left(
\{p,1-p\} \kron \{q,1-q\} \kron \{r,1-r\}
\right) \cdot 
(\{U_2\}_{ii})_{i=1}^8
\Big \}
\end{equation}
for $p,q,r\in[0,1]$.
\end{example}

\begin{figure}[ht!]
\centering
	\subfigure{
 \includegraphics[width=0.45\textwidth]{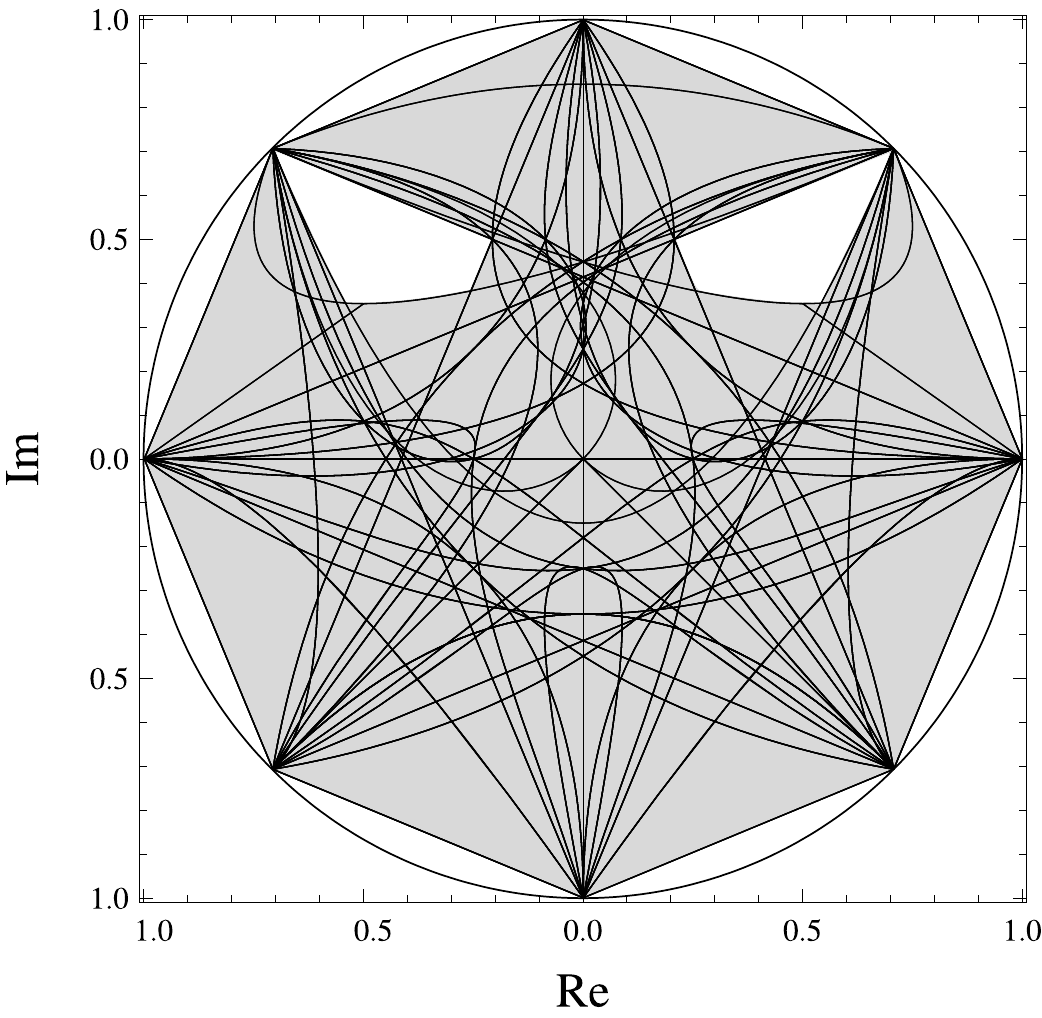}
 }
	\subfigure{
 \includegraphics[width=0.45\textwidth]{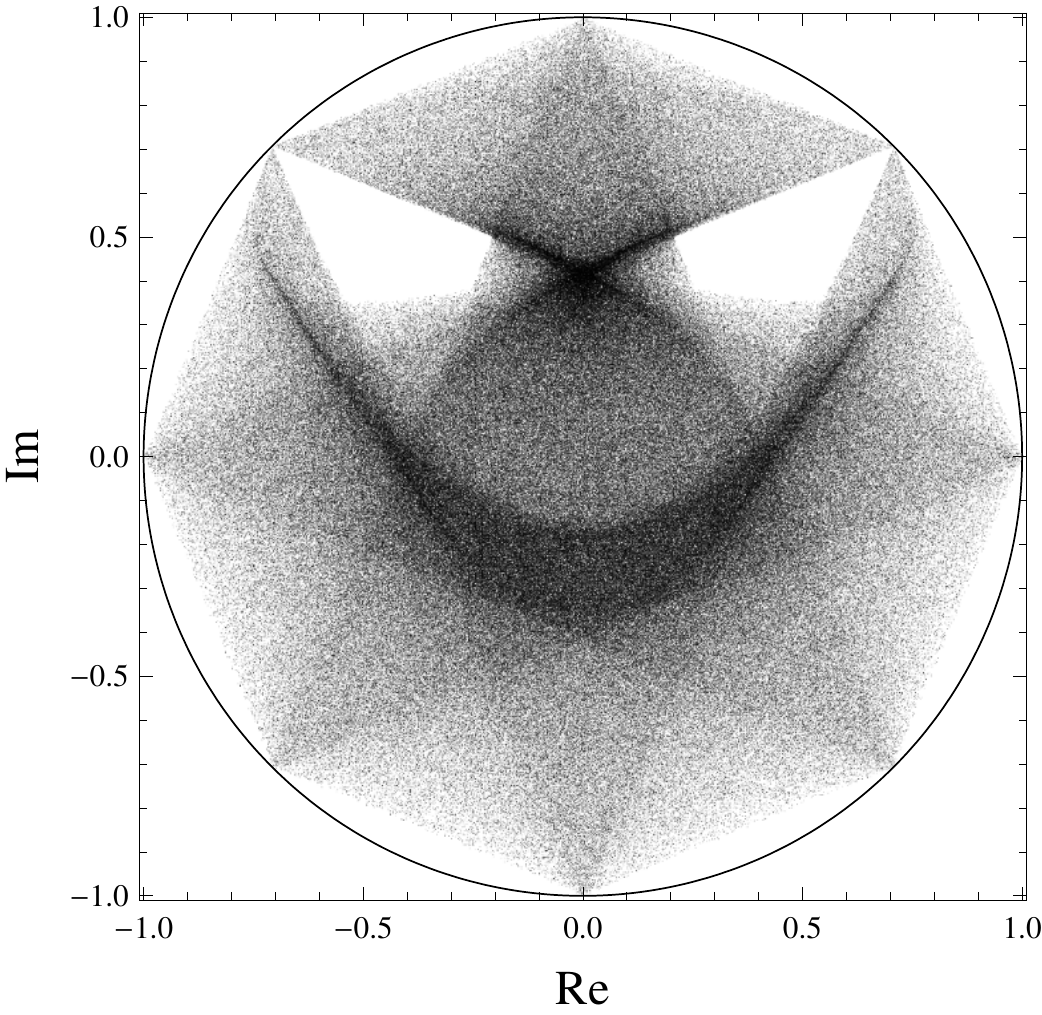}
 }
\caption{Product numerical range of matrix $A$ defined in (\ref{eqn:def-matrix-4}).
Panel a) shows an image of the product numerical range (plotted in gray) 
 obtained by parametrization~(\ref{eqn:par-4}),
while panel b) shows its image obtained by random sampling.
The standard numerical range forms a~regular octagon. 
}
\label{fig:genus2}
\end{figure}

\begin{example}[Genus 2]
Consider a~unitary operator of order $16$,
acting on a~four-qubit system. Assume that the 
corresponding matrix $A$ is diagonal in the product basis
$\{\ket{0000}, \dots , \ket{1111} \}$, 
\begin{equation}
\label{eqn:def-matrix-4}
A =\diag \bigl(
e^{\frac{i \pi }{4}},i,i,e^{\frac{3 i \pi }{4}},-1,e^{-\frac{3 i \pi }{4}},
e^{-\frac{i \pi }{4}},1,-1,e^{-\frac{i \pi }{4}},e^{-\frac{3 i \pi }{4}},1,
e^{\frac{3 i \pi }{4}},-i,-i,e^{\frac{i \pi }{4}} \bigr) .
\end{equation}

\rm\noindent Applying Proposition \ref{prop:diagonal-parametrization}, we parameterize
the product numerical range of this matrix in the following manner
\begin{equation}\label{eqn:par-4}
\ProductNumRange{A} =
\Big\{ z : z = 
\left(
\{p,1-p\} \kron \{q,1-q\} \kron \{r,1-r\} \kron \{s,1-s\}
\right) \cdot \lambda \Big \} 
 \end{equation}
for $ p,q,r,s\in[0,1]$, where '$\cdot$' denotes the scalar product.
Product numerical range of the matrix $A$ acting on the four--partite system 
forms a~set of genus $2$, as demonstrated in Figure~\ref{fig:genus2}. 
\end{example}

\section{Concluding Remarks}

In this work we established basic properties of product numerical range (PNR)
of operators acting on the Hilbert space with a tensor product structure.
Even though the definition of this extension of the standard numerical range
is rather straightforward, characterizing this set for a general operator
occurs to be a difficult problem. We have not managed therefore to present
its general solution, but we obtained concrete bounds for this set
(from inside and from outside) useful in various cases.
Several special cases of this general problem are relevant 
in view of practical applications of product numerical range
in the theory of quantum information \cite{Gxx10}.
To stimulate interest of the mathematical community in this subject 
we conclude the work by providing a list of open problems and sketching
their motivation stemming from physics.

\smallskip 

1. For a given Hermitian operator $H$ acting on a bi-partite Hilbert space,
   ${\cal H}_N={\cal H}_K \otimes {\cal H}_M$ find its local numerical range
   $[\lambda_{\rm min}^{\otimes}, \lambda_{\rm max}^{\otimes}]$.

\noindent Also precise bounds for PNR in this case would be useful. Of particular
   importance would be methods of establishing whether for a given $H$ the inequality
     $\lambda_{\rm min}^{\otimes} \ge 0$ is satisfied. This
     statement is equivalent of saying that the operator $H$ is {\it block positive}
     (with respect to the given decomposition of the space ${\cal H}_N$),
       hence the corresponding quantum map is positive \cite{Ja72,SZ09}.

\smallskip 

2. Find the local product range $\Lambda^{\otimes}(U)$ for a possibly large class of 
   unitary operators acting on a Hilbert space with the tensor product structure.

\noindent It is equally important to establish bounds applicable for PNR of an arbitrary  unitary $U$, suitable to answer the question if $0 \in \Lambda^{\otimes}(U)$.
This very issue occurs to be crucial in solving physical problems
related to local distinguishability of unitary quantum gates \cite{DFY08,Gxx10}.

\smallskip 

3. For a given operator $X$ acting on a Hilbert space with the tensor product structure
   find bounds for minimum of the modulus $|z|$, such that $z \in \Lambda^{\otimes}(X)$.

\noindent This question is related to finding the optimal gate fidelity \cite{Si09}
or optimizing local fidelity between two arbitrary quantum  states \cite{Gxx10}.

\smallskip 

4. Analyze the genus of PNR in the general case of operators acting on an $k$--partite Hilbert   space, ${\cal H}_1 \otimes \dots \otimes {\cal H}_k$.

\noindent In particular verify the conjecture that the genus $g$ is not larger than $k-2$.
       If this conjecture does not hold in general
 one may still try to check whether it is true for operators that are $k$-fold tensor products. This issue is then related to the question what the genus
 of the Minkowski product of $k$ convex sets on the complex plane is.

As a more specific question, one may ask the following:
\smallskip

5. Is $\Lambda^{\otimes}\left(A\otimes B\right)$ star-shaped for arbitrary operators $A$, $B$ acting on $\mathcal{H}_K$, $\mathcal{H}_M$, respectively?

We have already mentioned that problem in Section \ref{sec:tens_prod}.

\medskip

Furthermore, one may also consider similar problems in a more general set-up
 by studying corresponding product analogues of
other generalizations of the notion of numerical range.
For instance, mathematical results on 
`local (product) $C$--numerical range' and 
`local (product) $C$--numerical radius' 
\cite{dirr08relative,thomas08significance,SHGDH08},
  separable numerical range \cite{Gxx10}
and `Liouville numerical range' \cite{Si09}
will  find direct applications in several problems
in theoretical physics.

\section*{Acknowledgements}
It is a~pleasure to thank P.~Horodecki, C.K.~Li and T.~Schulte-Herbr{\"u}ggen
for fruitful discussions and to G.~Dirr, J. Gruska and  M.B. Ruskai
for helpful remarks. We acknowledge the financial support by
the Polish Ministry of Science and Higher Education under the grants number N519
012 31/1957 and DFG-SFB/38/2007, by the European research program COCOS and by
the Polish research network LFPPI. One of the authors (\L.S.) acknowledges that the project was operated within the Foundation
for Polish Science International Ph.D. Projects Program co-financed
by the European Regional Development Fund covering, under the agreement
no. MPD/2009/6, the Jagiellonian University International Ph.D. Studies in
Physics of Complex Systems.

\appendix
\section{Proof of Proposition \ref{productoftwo}}
\label{appendix_b}

 As a~consequence of Lemma \ref{lemmastar}, we can limit our proof to the 
case $0\not\in\LambdaRm\left(A_1\right)\cup\LambdaRm\left(A_2\right)$. Rotations of
 $\LambdaRm\left(A_1\right)$ or $\LambdaRm\left(A_2\right)$ in the complex plane cannot affect the genus of 
$\LambdaRm\left(A_1\otimes A_2\right)=\LambdaRm\left(A_1\right) \mprod \LambdaRm\left(A_2\right)$. 
Thus we may assume for the beginning that 
$\setR_+\cap\left(\LambdaRm\left(A_1\right)\cup\LambdaRm\left(A_2\right)\right)=\emptyset$.
We can achieve this by suitably rotating $\LambdaRm\left(A_1\right)$ and $\LambdaRm\left(A_2\right)$, because
 $\setR_+\cap\LambdaRm\left(A_j\right)=\emptyset$ or $-\setR_+\cap\LambdaRm\left(A_j\right)=\emptyset$
 for either $j$, as a~consequence of the convexity of $A_j$ and the assumption that $0\not\in A_j$. 
Now we can rotate $\LambdaRm\left(A_1\right)$ and $\LambdaRm\left(A_2\right)$ clockwise, so that $\phi=0$ becomes
 the minimal number $\phi$ in $\left[0,2\pi\right)$ for which $\setR_+e^{i\phi}\cap\LambdaRm\left(A_1\right)\neq\emptyset$, 
as well as the minimal number $\phi$ for which $\setR_+e^{i\phi}\cap\LambdaRm\left(A_2\right)\neq\emptyset$.
By our construction, keeping in mind that the sets $A_j$ are closed, there exist $\varepsilon_1,\varepsilon_2>0$
 such that $\setR_+e^{-i\varepsilon}\cap\LambdaRm\left(A_j\right)=\emptyset\forall_{\varepsilon\in\left(0,\varepsilon_j\right)}$
 for $j=1,2$. Using the convexity of $\LambdaRm\left(A_j\right)$ for $j=1,2$ again, we see that both $\LambdaRm\left(A_j\right)$ must 
be contained in the closed upper half-plane. Given $\setR_+\cap\LambdaRm\left(A_j\right)\neq\emptyset$, by the same argument as above we get 
$-\setR_+\cap\LambdaRm\left(A_j\right)=\setR_+e^{i\pi}\cap\LambdaRm\left(A_j\right)=\emptyset$. Closeness of $A_j$ implies that 
\begin{equation}
\label{phimax}
 \phi^j_+ =\max\left\{\phi\in\left[0,2\pi\right)|\setR_+e^{i\phi}\cap\LambdaRm\left(A_j\right)\neq\emptyset\right\}<\pi ,
\end{equation}
for $j=1,2$. Since $0\not\in \LambdaRm(A_1)\cup \LambdaRm(A_2)$ and $\phi^j_+<2\pi$, we can describe the sets $A_j$ in the 
log-polar coordinates $\Xi\left(z\right)=\left(\log|z|,\textnormal{Arg}\left(z\right)\right)$. Because the sets $A_j$ 
are convex, $\setR_+e^{i\phi}\cap\LambdaRm\left(A_j\right)$ must be a~closed nonempty interval or a~point for 
arbitrary $\phi$ and~$j$. Thus in the new coordinates we can write $\LambdaRm\left(A_j\right)$ as
\begin{equation}
\label{Ajribbon}
 \Xi \left(\LambdaRm\left(A_j\right)\right)=\bigcup_{\xi\in\left[0,\phi^j_+
\right]}\left[r^j_-\left(\xi\right),r^j_+\left(\xi\right)\right]\times\left\{\xi\right\}
\end{equation}
Using the convexity and closeness of $A_j$, it is easy to prove that the functions
 $r^j_{\pm}$ are continuous, which implies that \eqref{Ajribbon} is 
a simply connected subset of $\R^2$. Because the multiplication of complex numbers 
$z_1=|z_1|e^{i\phi_1},z_2=|z_2|e^{i\phi_2}$ corresponds to the addition of their arguments 
$\phi_1,\phi_2$ as well as the addition of the logarithms of their modules, $\log|z_1|$ and $\log|z_2|$, it is straightforward to
 write an explicit formula for $\LambdaRm\left(A_1\right)\LambdaRm\left(A_2\right)$ in the coordinates $\Xi$,
\begin{equation}
\label{productcoords}
 \bigcup_{\xi\in\left[0,\phi^1_++\phi^2_+\right]}\bigcup_{\varphi\in\left[\max\left(0,\xi-\phi^2_+
\right),\min\left(\phi^1_+,\xi\right)\right]}
\left[r^1_+\left(\varphi\right)r^2_+\left(\xi-\varphi\right),r^1_-\left(\varphi\right)r^2_-\left(\xi-\varphi\right)\right]
\times\left\{\xi\right\}
\end{equation}
or
\begin{equation}
\label{productcoordswithR}
 \bigcup_{\xi\in\left[0,\phi^1_++\phi^2_+\right]}
\left[R_+\left(\xi\right),R_-\left(\xi\right)\right]
\times\left\{\xi\right\},
\end{equation}
where the functions $R^{\pm}$ are defined in the following way, \begin{eqnarray}\label{defofRs}R^+\left(\xi\right)&=\max_{\left[\max\left(0,\xi-\phi^2_+
\right),\min\left(\phi^1_+,\xi\right)\right]}\left(r^1_+\left(\varphi\right)r^2_+\left(\xi-\varphi\right)\right),\\R^-\left(\xi\right)&=\min_{\left[\max\left(0,\xi-\phi^2_+
\right),\min\left(\phi^1_+,\xi\right)\right]}\left(r^1_-\left(\varphi\right)r^2_-\left(\xi-\varphi\right)\right).\end{eqnarray}

It turns out that $R_+$ and $R_-$ are continuous and the interval $\left[R_-\left(\xi\right),R_+\left(\xi\right)\right]$ is nonempty for arbitrary $\xi\in\left[0,\phi^1_++\phi^2_+\right]$. Therefore \eqref{productcoords} is a~simply connected subset of $\setR\times\left[0,\phi^1_++\phi^2_+\right]$. 
Since $\phi^1_++\phi^2_+<2\pi$, $\Xi$ is a~homomorphism between $\setR\times\left[0,\phi^1_++\phi^2_+\right]$ 
and $\Xi^{-1}\left(\setR\times\left[0,\phi^1_++\phi^2_+\right]\right)$. 
Therefore $\LambdaRm\left(A_1\right) \mprod \LambdaRm\left(A_2\right)$ is simply connected, as a~preimage 
of a~simply connected set by the isomorphism $\Xi$. This is what we wanted to prove.
\section{Proof of Proposition \ref {prop:diagonal-parametrization}}
\label{sec:proof-th-diagonal-parametrization}
Each product vector $\ket{x} = \ket{x_1} \kron \ket{x_2} \dots \kron \ket{x_k}$ of unit norm can be rewritten in its computational basis,
\begin{equation}
(x_{1,0} \ket{0}+ \dots + x_{1,m_1 -1} \ket{m_1-1}) \kron \dots \kron (x_{k,0} \ket{0}+ \dots + x_{k,m_k -1} \ket{m_k-1}),
\end{equation}
where 
\begin{equation}\label{eqn:warunekNaSumeKwadratow}
|x_{r,0}|^2 + \dots + |x_{r,m_{r}-1}|^2 = 1 \; \; \text{ for } r=1 \dots k.
\end{equation}
Thus we have
\begin{eqnarray}
\nonumber
 \bra{x} A \ket{x} &=& 
\sum_{l_1=0}^{m_1-1} \sum_{l_2=0}^{m_2-1} \dots \sum_{l_k=0}^{m_k-1}
\sum_{s_1=0}^{m_1-1} \sum_{s_2=0}^{m_2-1} \dots \sum_{s_k=0}^{m_k-1} \\
&&
	\bar{x}_{1,l_1}\bar{x}_{2,l_2} \dots \bar{x}_{k,l_k} 
	x_{1,s_1} x_{2,s_2} \dots x_{k,s_k} 
	\bra{l_1 l_2 \dots l_k} A \ket{s_1 s_2 \dots s_k}.
\end{eqnarray}
Now we must note that, since $A$ is diagonal with respect to the product computational basis, we have
\begin{equation}
\bra{l_1 l_2 \dots l_k} A \ket{s_1 s_2 \dots s_k} = \lambda_{l_1 l_2 \dots l_k} \delta_{\{ l_1 l_2 \dots l_k , s_1 s_2 \dots s_k \}}.
\end{equation}
Thus we get
\begin{equation}
 \bra{x} A \ket{x} = 
\sum_{l_1=0}^{m_1-1} \sum_{l_2=0}^{m_2-1} \dots \sum_{l_k=0}^{m_k-1}
	|x_{1,l_1}|^2 |x_{2,l_2}|^2 \dots |x_{k,l_k}|^2
	\lambda_{l_1 l_2 \dots l_k}.
\end{equation}
If we denote $p^{(r)}_{i} := |x_{r,i}|^2$, then we have $p^{(r)}_0 + p^{(r)}_1 + \dots + p^{(r)}_{m_r-1}= 1$ and
\begin{equation}
 \bra{x} A \ket{x} = 
\sum_{l_1=0}^{m_1-1} \sum_{l_2=0}^{m_2-1} \dots \sum_{l_k=0}^{m_k-1}
	p^{(1)}_{l_1} p^{(2)}_{l_2} \dots p^{(k)}_{l_k}
	\lambda_{l_1 l_2 \dots l_k}.
\end{equation}
Now it is easy to notice that if we take all possible product states, 
we will obtain all possible decompositions 
\begin{equation}
 p^{(r)}_0,p^{(r)}_1,\dots,p^{(r)}_{m_r} \geq 0, \; 
p^{(r)}_0 + p^{(r)}_1 + \dots + p^{(r)}_{m_r}= 1
\end{equation}
for $r = 1, \dots k$.
Thus the proof is complete.
\halmos

\end{document}